\documentclass[12pt]{article}

\usepackage{amsmath, epsfig, cite}
\usepackage{amssymb}
\usepackage{amsfonts}
\usepackage{latexsym}
\usepackage{amsthm}

\makeatletter
\renewcommand{\@seccntformat}[1]{{\csname the#1\endcsname}.\hspace{.5em}}
\makeatother

\newtheorem{thm}{Theorem}[section]

\newtheorem{cor}[thm]{Corollary}
\newtheorem{conj}[thm]{Conjecture}
\newtheorem{lem}[thm]{Lemma}

\newcommand{\pf}{\noindent{\it Proof.} }

\renewcommand{\qed}{\hfill$\Box$\medskip}

\numberwithin{equation}{section}

\begin{document}

\renewcommand{\thefootnote}{*}

\begin{center}
{\Large\bf $q$-Analogues of two ``divergent" Ramanujan-type\\[5pt] supercongruences}
\end{center}

\vskip 2mm \centerline{Victor J. W. Guo  }
\begin{center}
{\footnotesize School of Mathematical Sciences, Huaiyin Normal University, Huai'an, Jiangsu 223300,\\
 People's Republic of China\\
{\tt jwguo@hytc.edu.cn } }
\end{center}

\vskip 0.7cm \noindent{\small{\bf Abstract.}} Guillera and Zudilin proved three ``divergent" Ramanujan-type
supercongruences by means of the Wilf-Zeilberger algorithmic technique. In this paper, we prove $q$-analogues of two of them
via the $q$-WZ method. Additionally, we give $q$-analogues of two related congruence of Sun; one is confirmed and the other is conjectural.



\section{Introduction}
By using the Wilf--Zeilberger method, Guillera and Zudilin \cite{GuZu} proved the following three supercongruences:
\begin{align}
\sum_{k=0}^{\frac{p-1}{2}} \frac{(\frac{1}{2})_k^3}{k!^3}(3k+1)2^{2k} &\equiv p\pmod{p^3}\quad\text{for $p>2$}, \label{eq:div-1}\\[5pt]
\sum_{k=0}^{\frac{p-1}{2}} \frac{(\frac{1}{2})_k^5}{k!^5}(10k^2+6k+1)2^{2k} &\equiv p^2\pmod{p^5}\quad\text{for $p>3$}, \label{eq:div-2}\\[5pt]
\sum_{k=0}^{\frac{p-1}{2}} \frac{(\frac{1}{2})_k^3}{k!^3}(3k+1)(-1)^k 2^{3k} &\equiv p(-1)^{\frac{p-1}{2}}\pmod{p^3}\quad\text{for $p>2$}.  \label{eq:div-3}
\end{align}
Here, and throughout the paper, the letter $p$ always denotes a prime, and the Pochhammer symbol $(a)_b$ is used for denoting $\Gamma(a+b)/\Gamma(b)$ also in
the cases when $b$ is not a non-negative integer.
In the sprit of \cite{Zudilin}, the supercongruences \eqref{eq:div-1}--\eqref{eq:div-3} correspond to divergent Ramanujan-type series for $1/\pi$ or $1/\pi^2$, such as
\begin{align}
\sum_{k=0}^{\infty} \frac{(\frac{1}{2})_k^3}{(1)_k^3}(3k+1)2^{2k}\ \text{``="}\ \frac{-2i}{\pi},\qquad
\sum_{k=0}^{\infty} \frac{(\frac{1}{2})_k^3}{(1)_k^3}(3k+1)(-1)^k 2^{3k}\ \text{``="}\ \frac{1}{\pi}  \label{eq:div-4}
\end{align}
(see \cite[(47)]{GuZu}). The summations in \eqref{eq:div-4} have to be understood as the analytic continuation of the corresponding hypergeometric series. For instance,
the second formula in \eqref{eq:div-4} can be written as
\begin{align*}
\frac{1}{2\pi i}\int_{-i\infty}^{i\infty}\frac{(\frac{1}{2})_s^3}{(1)_s^3}\Gamma(-s)(3s+1) 2^{3s}ds=\frac{1}{\pi}.
\end{align*}
It is worth mentioning that Guillera \cite{Guillera3} has given proofs of several divergent hypergeometric formulas for $1/\pi$ and $1/\pi^2$ by using a version of the Wilf--Zeilberger method.

In a previous paper, motivated by Zudilin's work \cite{Zudilin},
the first author \cite{Guo2018} utilized the $q$-WZ method \cite{Mohammed,Zeilberger} to give the following $q$-analogue of a Ramanujan-type supercongruence of van Hamme \cite{Hamme}:
\begin{align}
\sum_{k=0}^{\frac{p-1}{2}}(-1)^k q^{k^2}[4k+1]\frac{(q;q^2)_k^3}{(q^2;q^2)_k^3}
\equiv [p]q^{\frac{(p-1)^2}{4}} (-1)^{\frac{p-1}{2}}\pmod{[p]^3}\quad\text{for $p>2$},  \label{eq:q-Hamme}
\end{align}
where the {\it $q$-shifted factorial} is defined by $(a;q)_0=1$ and $(a;q)_n=(1-a)(1-aq)\cdots (1-aq^{n-1})$ for $n\geqslant 1$,
while the {\it $q$-integer} is defined as $[n]=1+q+\cdots+q^{n-1}$ (see \cite{GR}).
We point out that, in this paper, two rational functions in $q$ are congruent
modulo a polynomial $P(q)$ if the numerator of their difference is congruent to $0$ modulo $P(q)$
in the polynomial ring $\mathbb{Z}[q]$ while the denominator is relatively prime to $P(q)$.
Note that some other interesting $q$-congruences can be found in \cite{Andrews,CP,GZ15,LPZ,SP,Tauraso1,Tauraso2}.

Although supercongruences have been widely studied by many mathematicians including
Beukers \cite{Beukers1}, Long and Ramakrishna \cite{LR}, Rodriguez-Villegas \cite{RV}, Z.-H. Sun \cite{SunZH}, Z.-W. Sun \cite{Sun4}, van Hamme \cite{Hamme}, and Zudilin \cite{Zudilin}, etc.,
there are still many problems on $q$-congruences which are worthwhile to investigate. In fact, it is not always easy to give $q$-analogues of known congruences,
and many congruences might have no $q$-analogues. Sometimes in order to prove an ordinary congruence (without $q$) we need to
establish its $q$-analogue \cite{CP}. Proving $q$-congruences requires a variety of methods, including properties of root of unity \cite{SP}, basic hypergeometric series identities
\cite{GZ15,LPZ,Tauraso1}, the $q$-Wilf-Zeilberger ($q$-WZ) method \cite{Guo2018,Tauraso2}. In many cases the method to prove an ordinary congruence cannot be generalized directly to prove
the corresponding $q$-congruence. For instance, no one knows how to extend the $p$-adic analysis in \cite{KLMSY} to the $q$-case.
On the other hand, sometimes studying $q$-congruences will enable us to find
new basic hypergeometric series identities \cite{GGZ,GZ15}. Moreover, combinatorial proofs of ordinary congruences are little known.
Finding out the $q$-analogues will help us to have a better understanding of the original congruences, and the theory
of integer partitions might be applied to give the desirable combinatorial interpretations.

This paper is a continuation of \cite{Guo2018}. We want to give $q$-analogues of \eqref{eq:div-1} and \eqref{eq:div-3} by the $q$-WZ method.
Recall that the $n$-th {\it cyclotomic polynomial} that $\Phi_n(q)$ may be defined as
\begin{align}
\Phi_n(q):=\prod_{\substack{1\leqslant k\leqslant n\\ \gcd(n,k)=1}}(q-\zeta^k),  \label{eq:phi}
\end{align}
where $\zeta$ is a $n$-th primitive root of unity. It is well known that $\Phi_n(q)$ is an irreducible polynomial with integer coefficients and
$\Phi_p(q)=[p]$.

Note that we can sum in \eqref{eq:div-1}--\eqref{eq:div-3} up to $p-1$, since the $p$-adic order of $(\frac{1}{2})_k/k!$ is $1$ for $\frac{p+1}{2}\leqslant k\leqslant p-1$.
Our $q$-analogues of \eqref{eq:div-1} can be stated as follows.

\begin{thm}\label{thm:q-Zudilin-0}
Let $n$ be a positive odd integer. Then
\begin{align}
\sum_{k=0}^{\frac{n-1}{2}}[3k+1]\frac{(q;q^2)_k^3 q^{-{k+1\choose 2} } }{(q;q)_k^2 (q^2;q^2)_k}
&\equiv [n]q^{\frac{1-n}{2}}\pmod{[n]\Phi_n(q)^2}, \label{eq:q-div-WZ-1} \\[5pt]
\sum_{k=0}^{n-1}[3k+1]\frac{(q;q^2)_k^3 q^{-{k+1\choose 2} } }{(q;q)_k^2 (q^2;q^2)_k}
&\equiv [n]q^{\frac{1-n}{2}}\pmod{[n]\Phi_n(q)^2}.  \label{eq:q-div-WZ-2}
\end{align}
\end{thm}

It is easy to see that, when $n=p$ is an odd prime, the congruences \eqref{eq:q-div-WZ-1} and \eqref{eq:q-div-WZ-2} are equivalent to each other,
since $\frac{(q;q^2)_k^3 }{(q;q)_k^2 (q^2;q^2)_k}\equiv 0\pmod{[p]}$ for
$\frac{p+1}{2}\leqslant k\leqslant p-1$. But for general $n$ they are clearly not equivalent.

Let $n=p^r$ be an odd prime power and $q=1$ in Theorem \ref{thm:q-Zudilin-0}. Noticing that $\Phi_{p^r}(1)=p$ and the denominator of the reduced form of
$\frac{(\frac{1}{2})_k^3}{k!^3}(3k+1)2^{2k}$ is relatively prime to $p$,
we obtain the following generalization of \eqref{eq:div-1}.
\begin{cor}
Let $p>2$ and $r>0$. Then
\begin{align}
\sum_{k=0}^{\frac{p^r-1}{2}} \frac{(\frac{1}{2})_k^3}{k!^3}(3k+1)2^{2k} &\equiv p^r \pmod{p^{r+2}}, \label{eq:div-gen-r-1} \\[5pt]
\sum_{k=0}^{p^r-1} \frac{(\frac{1}{2})_k^3}{k!^3}(3k+1)2^{2k} &\equiv p^r \pmod{p^{r+2}}.  \label{eq:div-gen-r-2}
\end{align}
\end{cor}

We also have two different $q$-analogues of \eqref{eq:div-3} as follows.
\begin{thm}\label{thm:q-Zudilin-3}
Let $n$ be a positive odd integer. Then
\begin{align}
\sum_{k=0}^{n-1}(-1)^k [3k+1]\frac{(q;q^2)_k^3}{(q;q)_k^3}
\equiv [n]q^{\frac{(n-1)^2}{4}}(-1)^{\frac{n-1}{2}} \pmod{[n]\Phi_n(q)^2}, \label{eq:q-Zudilin-3}
\end{align}
\end{thm}
\begin{conj}\label{conj:q-Zudilin-3}
Let $n$ be a positive odd integer. Then
\begin{align}
\sum_{k=0}^{\frac{n-1}{2}}(-1)^k [3k+1]\frac{(q;q^2)_k^3}{(q;q)_k^3}
\equiv [n]q^{\frac{(n-1)^2}{4}}(-1)^{\frac{n-1}{2}} \pmod{[n]\Phi_n(q)^2}. \label{eq:conj-Zudilin-3}
\end{align}
\end{conj}

It should be mentioned that Sun \cite[Conjecture 5.1(ii)]{Sun4} discovered refinements of \eqref{eq:div-1} and \eqref{eq:div-3} modulo $p^4$,
which has been recently confirmed by Mao and Zhang \cite{MZ} and Chen, Xie, and He\cite{CXH}, respectively. Moreover, Sun \cite[Conjecture 5.1(i)]{Sun4} also proposed the following conjecture:
\begin{align}
\sum_{k=0}^n (3k+1){2k\choose k}^3 16^{n-k}    &\equiv 0 \mod{4(2n+1){2n\choose n}}, \label{eq:sun-bino-1}\\[5pt]
\sum_{k=0}^n (3k+1){2k\choose k}^3 (-8)^{n-k} &\equiv 0 \mod{4(2n+1){2n\choose n}}, \label{eq:sun-bino-2}
\end{align}
which has been recently proved by Mao and Zhang \cite{MZ} and He \cite{He}, respectively. In this paper, we shall give a $q$-analogue of \eqref{eq:sun-bino-1} as follows.

\begin{thm}\label{thm:q-Zudilin-MZ}
Let $n$ be a positive integer. Then
\begin{align}
\sum_{k=0}^{n}(-1)^k [3k+1]{2k\brack k}^3 \frac{(-q;q)_n^3}{(-q;q)_k^3}
\equiv 0\pmod{(1+q^n)^2[2n+1]{2n\brack n}}, \label{eq:q-Sun-2}
\end{align}
where the {\it $q$-binomial coefficients} ${M\brack N}$
are defined by
$$
{M\brack N}={M\brack N}_q
=\begin{cases}\displaystyle\frac{(q;q)_M}{(q;q)_N(q;q)_{M-N}} &\text{if $0\leqslant N\leqslant M$,} \\[5pt]
0 &\text{otherwise.}
\end{cases}
$$
\end{thm}

From \eqref{eq:q-Sun-2}, we can easily deduce the following generalization of \eqref{eq:div-3}.
\begin{cor}Let $p>2$ and $r>0$. Then
\begin{align}
\sum_{k=0}^{\frac{p^r-1}{2}} \frac{(\frac{1}{2})_k^3}{k!^3}(3k+1)(-1)^k 2^{3k} &\equiv p^r(-1)^{\frac{p-1}{2}}\pmod{p^{r+2}}\quad\text{for $p>2$}.  \label{eq:div-3-pr}
\end{align}
\end{cor}

We also have the following conjecture on a $q$-analogue of \eqref{eq:sun-bino-2}.

\begin{conj}\label{conj:q-Zudilin-MZ}
Let $n$ be a positive integer. Then
\begin{align}
\sum_{k=0}^{n}[3k+1]{2k\brack k}^3 \frac{(-q;q)_n^4}{(-q;q)_k^4} q^{-{k+1\choose 2}}
\equiv 0\pmod{(1+q^n)^2[2n+1]{2n\brack n}}.\label{eq:q-Sun-1}
\end{align}
\end{conj}

Sun and Tauraso (a special case of \cite[Theorem 1.3]{ST}) proved that
\begin{align}
\sum_{k=1}^{p-1}{2k\choose k}\frac{1}{k}\equiv 0\pmod{p^2}\quad\text{for $p>3$}.  \label{eq:st}
\end{align}
In this paper, we shall give the following $q$-analogue of \eqref{eq:st}.
\begin{thm}\label{thm:q-st}
Let $n$ be a positive odd integer. Then
\begin{align}
\sum_{k=1}^{n-1}\frac{[3k]}{[2k]^2}{2k\brack k}q^{-{k\choose 2}}
\equiv \frac{(n^2-1)(1-q)^2}{24}[n] \pmod{\Phi_n(q)^2}.  \label{eq:q-st}
\end{align}
\end{thm}

It is easy to see that when $n=p>3$ is a prime and $q=1$ the congruence \eqref{eq:q-st} reduces to \eqref{eq:st}.
However, we cannot obtain any interesting congruence from \eqref{eq:q-st} for the case $q=1$ and $n=p^r$ with $r>1$, since when $n=p^r$
the denominator of the reduced form of the left-hand side of
\eqref{eq:q-st} has factors $\Phi_{p^j}(q)$ for $1\leqslant j\leqslant r-1$.

The rest of the paper is organized as follows. In Section 2 we present auxiliary congruences, some of which are of interest on their own.
We shall prove Theorems \ref{thm:q-Zudilin-0} and \ref{thm:q-Zudilin-3} in Sections 3 and 4 respectively using the $q$-WZ method.
In Section 5, we shall prove Theorem \ref{thm:q-Zudilin-MZ} using the same $q$-WZ pair in Section 4. We give a short proof of Theorem~\ref{thm:q-st} in Section 6.
The final section, Section 7, provides some remarks on Conjecture \ref{conj:q-Zudilin-MZ} and proposes more conjectures for further study, including a refinement of the congruence \eqref{eq:q-div-WZ-2}.

\section{Pecongruences }
In this section, we summarize our needs for proving the congruences of Theorems \ref{thm:q-Zudilin-0} and \ref{thm:q-Zudilin-MZ}.
Recall that Staver's identity \cite{Staver} can be written as
\begin{align*}
\sum_{k=1}^{N}{2k\choose k}\frac{1}{k}=\frac{N+1}{3}{2N+1\choose N}\sum_{k=1}^{N}\frac{1}{k^2{N\choose k}^2}. 
\end{align*}
It plays an important part in Zudilin's proof of \eqref{eq:div-1}, and has also been utilized by Sun \cite{Sun4}, Sun and Tauraso \cite{ST}, and Mao and Zhang \cite{MZ}
to prove certain supercongruences. The following result is a $q$-analogue of Staver's identity.

\begin{lem}Let $n$ be a positive integer. Then
\begin{align}
\sum_{k=1}^{n}\frac{[3k]}{[2k]^2}{2k\brack k}q^{-{k\choose 2}}
=[n+1]{2n+1\brack n}\sum_{k=1}^{n}\frac{q^{-{n-2k+1\choose 2}}}{[2k]^2{n\brack k}_{q^2}^2}. \label{eq:q-Staver}
\end{align}
\end{lem}
\pf We use the same technique  as Paule \cite{Paule}, who introduced symmetry factors to simplify the proofs of terminating $q$-hypergeometric identities.
It is easy to see that
$$
\sum_{k=1}^{n}\frac{(1-q^{2k-n-1})q^{-{n-2k+1\choose 2}}}{[2k]^2{n\brack k}_{q^2}^2}=0,
$$
since the $k$-th summand and the $n+1-k$ summand cancel each other. Hence, the identity \eqref{eq:q-Staver} is equivalent to
\begin{align}
\sum_{k=1}^{n}\frac{[3k]}{[2k]^2}{2k\brack k}q^{-{k\choose 2}}
=\frac{[n+1]}{2}{2n+1\brack n}\sum_{k=1}^{n}\frac{(1+q^{2k-n-1})q^{-{n-2k+1\choose 2}}}{[2k]^2{n\brack k}_{q^2}^2}.
\label{eq:new-q-Staver}
\end{align}

Let
\begin{align*}
F(n,k)&=\frac{[n+1]}{2}{2n+1\brack n}\frac{(1+q^{2k-n-1})q^{-{n-2k+1\choose 2}}}{[2k]^2{n\brack k}_{q^2}^2},\\[5pt]
G(n,k)&=-\frac{[3n-2k+5]}{2}{2n+1\brack n}\frac{(1+q^{n+1})q^{-{n-2k+3\choose 2}}}{[2k]^2{n+1\brack k}_{q^2}^2}.
\end{align*}
Then we can check that
\begin{align}
F(n+1,k)-F(n,k)
=G(n,k+1)-G(n,k).  \label{eq:diff-fnk-gnk}
\end{align}
Namely, the functions $F(n,k)$ and $G(n,k)$ form a $q$-WZ pair. Summing \eqref{eq:diff-fnk-gnk} over $k$ from $1$ to $n$, we obtain
\begin{align*}
&\hskip -3mm \sum_{k=1}^{n+1}F(n+1,k)-\sum_{k=1}^{n}F(n,k)  \\[5pt]
&=F(n+1,n+1)+G(n,n+1)-G(n,1)  \\[5pt]
&=\frac{[n+2]{2n+3\brack n+1}(1+q^n)q^{-{-n\choose 2}}}{2[2n+2]^2}-\frac{[n+3]{2n+1\brack n}(1+q^{n+1})q^{-{1-n\choose 2}}}{2[2n+2]^2} \\[5pt]
&\quad{}+\frac{[3n+3]{2n+1\brack n}(1+q^{n+1})q^{-{n+1\choose 2}}}{2[2n+2]^2} \\[5pt]
&=\frac{[3n+3]}{[2n+2]^2}{2n+2\brack n+1}q^{-{n+1\choose 2}}.
\end{align*}
It is clear that the identity \eqref{eq:new-q-Staver} immediately follows  from the above recurrence by induction on $n$.
\qed

We now give a $q$-analogue of \cite[(13)]{GuZu}.
\begin{lem}Let $n$ be a positive odd integer. Then
\begin{align}
\sum_{k=1}^{\frac{n-1}{2}}\frac{[3k]}{[2k]^2}{2k\brack k}q^{-{k\choose 2}}\equiv 0 \pmod{\Phi_n(q)},  \label{eq:3k2kk}
\end{align}
\end{lem}
\pf Replacing $n$ by $\frac{n-1}{2}$ in \eqref{eq:q-Staver} and noticing that
${n\brack \frac{n-1}{2}}\equiv 0\pmod{\Phi_n(q)}$ and $\gcd([2k]{\frac{n-1}{2}\brack k}_{q^2},\Phi_n(q))=1$ for $1\leqslant k\leqslant\frac{n-1}{2}$, we obtain \eqref{eq:3k2kk}.
\qed

We need two auxiliary lemmas on properties of $q$-factorials.
\begin{lem}\label{eq:lem-2nn}
Let $n$ be a positive odd integer, and let $k$ be a non-negative integer. Then
\begin{align}
\frac{(q;q^2)_{n}(q^{2k+1};q^2)_{n-1}^2}{(q;q)_{n-1}^3} &\equiv 0\pmod{[n]\Phi_n(q)^2}. \label{eq:mod-q2n}
\end{align}
Moreover, if $k\leqslant\frac{n-1}{2}$, then
\begin{align}
\frac{(q;q^2)_{(n+1)/2}(q^{2k+1};q^2)_{(n-1)/2}^2}{(q;q)_{(n-1)/2}^3} &\equiv 0\pmod{[n]\Phi_n(q)^2}. \label{eq:mod-q2n-3}
\end{align}
\end{lem}
\pf It is well known that
\begin{align}
q^n-1 =\prod_{d|n}\Phi_d(q).
\end{align}
Therefore,
\begin{align*}
(q;q)_{n-1}&=(-1)^{n-1}\prod_{d=1}^{n-1}\Phi_{d}(q)^{\lfloor\frac{n-1}{d}\rfloor},  \\[5pt]
(q;q^2)_n
&=\frac{(q;q)_{2n}}{(q^2;q^2)_n}
=(-1)^{n} \prod_{d=1}^{n} \Phi_{2d-1}(q)^{\lfloor\frac{2n}{2d-1}\rfloor-\lfloor\frac{n}{2d-1}\rfloor},   \\[5pt]
(q^{2k+1};q^2)_n
&=\frac{(q;q^2)_{n+k}}{(q;q^2)_k}=(-1)^{n} \prod_{d=1}^{n+k}
\Phi_{2d-1}(q)^{\lfloor\frac{2n+2k}{2d-1}\rfloor+\lfloor\frac{k}{2d-1}\rfloor-\lfloor\frac{n+k}{2d-1}\rfloor-\lfloor\frac{2k}{2d-1}\rfloor },
\end{align*}
where $\lfloor x\rfloor$ denotes the greatest integer less than or equal to $x$.

We now suppose that $t|n$ and $t>1$. Since $n$ is odd, we know that $t$ is also odd. Hence, the exponent of $\Phi_t(q)$ in
$\frac{(q^{2k+1};q^2)_n}{(q;q)_{n-1}}$ is
$\lfloor\frac{2n+2k}{t}\rfloor+\lfloor\frac{k}{t}\rfloor-\lfloor\frac{n+k}{t}\rfloor
-\lfloor\frac{2k}{t}\rfloor-\lfloor\frac{n-1}{t}\rfloor=\frac{n}{t}-\lfloor\frac{n-1}{t}\rfloor=1$.
Since $[n]=\prod_{t|n,\ t>1}\Phi_t(q)$, we obtain
\begin{align}
\frac{(q^{2k+1};q^2)_{n}}{(q;q)_{n-1}} &\equiv 0\pmod{[n]}. \label{eq:mod-q2n-1}
\end{align}
It is easy to see that
\begin{align}
\frac{(q^{2k+1};q^2)_{n-1}}{(q;q)_{n-1}} &\equiv 0\pmod{\Phi_n(q)}. \label{eq:mod-q2n-2}
\end{align}
Moreover, using the following inequality
\begin{align}
 \lfloor 2x+2y\rfloor+\lfloor y\rfloor\geqslant \lfloor x\rfloor+\lfloor x+y\rfloor+\lfloor 2y\rfloor \label{eq:ineq}
\end{align}
(see \cite{Bober,WZ} for more such inequalities and related results), we can show that
the denominator of the reduced form of \eqref{eq:mod-q2n-2} is a product of even-th cyclotomic polynomials,
and is therefore relatively prime to $[n]$.
The proof of \eqref{eq:mod-q2n} then follows from \eqref{eq:mod-q2n-1} with $k=0$ and \eqref{eq:mod-q2n-2}.

Similarly, for positive odd integer $t$, the exponent of $\Phi_t(q)$ in
$\frac{(q^{2k+1};q^2)_{(n-1)/2}}{(q;q)_{(n-1)/2}}$ is
$\lfloor\frac{n-1+2k}{t}\rfloor+\lfloor\frac{k}{t}\rfloor-\lfloor\frac{n-1+2k}{2t}\rfloor
-\lfloor\frac{2k}{t}\rfloor-\lfloor\frac{n-1}{2t}\rfloor\geqslant 0$ by \eqref{eq:ineq}. Furthermore, if $k\leqslant\frac{n-1}{2}$, then
the exponent of $\Phi_n(q)$ in
$\frac{(q^{2k+1};q^2)_{(n-1)/2}}{(q;q)_{(n-1)/2}}$ is
$\lfloor\frac{n-1+2k}{n}\rfloor+\lfloor\frac{k}{n}\rfloor-\lfloor\frac{n-1+2k}{2n}\rfloor
-\lfloor\frac{2k}{n}\rfloor-\lfloor\frac{n-1}{2n}\rfloor=1$.
The proof of \eqref{eq:mod-q2n-3} then follows from $(q;q^2)_{(n+1)/2}=(1-q)[n](q;q^2)_{(n-1)/2}$.
\qed

\begin{lem}\label{lem:big-qbino}
Let $n$ be a positive odd integer. Then for $k=1,\ldots,\frac{n-1}{2}$ we have
\begin{align}
[n]{2n-2k\brack n-1}\frac{(q;q^2)_{n}(q;q^2)_{n-k} }{(q;q)_{n} (q^2;q^2)_{n-k}} \equiv 0\pmod{[n]\Phi_{n}(q)^2}.  \label{eq:big-qbino}
\end{align}
\end{lem}
\pf We can write the left-hand side of \eqref{eq:big-qbino} as
\begin{align*}
[n]{2n-2k\brack n-1}{2n\brack n}{2n-2k\brack n-k}\frac{1}{(-q;q)_n (-q;q)_{n-k}^2}.
\end{align*}
For $k=1,\ldots,\frac{n-1}{2}$, we have $2n-2k>n$ and $n-2k+1<n$. Therefore, by the definition of $q$-binomial coefficients, we immediately get
$$
{2n-2k\brack n-1}\equiv {2n-2k\brack n-k}\equiv 0\pmod{\Phi_n(q)}.
$$
The proof then follows from the fact that $[n]$ is relatively prime to $(-q;q)_n (-q;q)_{n-k}^2$.
\qed

We need the following result for the proof of Theorem \ref{thm:q-st}.
\begin{lem}Let $n$ be a positive integer. Then
\begin{align}
\sum_{k=1}^{n-1}\frac{q^{k}}{[2k]^2}\equiv \frac{(n^2-1)(1-q)^2}{24} \pmod{\Phi_n(q)}.  \label{eq:sp-sim}
\end{align}
\end{lem}
\pf The proof is similar to that of \cite[Lemma 2]{SP}.  For the sake of completeness,
we provide it here. Since $\frac{1}{[2k]}=\frac{1-q}{1-q^{2k}}$, the congruence \eqref{eq:sp-sim} is equivalent to
\begin{align}
G(q):=\sum_{k=1}^{n-1}\frac{q^{k}}{(1-q^{2k})^2}\equiv \frac{n^2-1}{24} \pmod{\Phi_n(q)}.  \label{eq:sp-sim}
\end{align}
Let $\zeta=e^{\frac{2\pi i}{n}}$ be a $n$-th primitive root of unity. By \eqref{eq:phi}, to prove \eqref{eq:sp-sim}, it suffices to show that
\begin{align*}
G(\zeta^m)=\frac{n^2-1}{24}
\end{align*}
for all positive integers $m\leqslant n-1$ such that $\gcd(m,n)=1$. Since $n$ is odd, it is easy to see that
\begin{align*}
G(\zeta^m)=\sum_{k=1}^{n-1}\frac{\zeta^{mk}}{(1-\zeta^{2mk})^2}= \sum_{k=1}^{n-1}\frac{\zeta^{k}}{(1-\zeta^{2k})^2}=G(\zeta),
\end{align*}
provided that $\gcd(m,n)=1$.

We now define
\begin{align*}
G(q,z):=\sum_{k=1}^{n-1}\frac{q^{k}}{(1-q^{2k}z)^2}.
\end{align*}
As
$$
\sum_{k=1}^{n-1}\zeta^{jk}
=\begin{cases}
n-1,&\text{if $n\mid j$},\\
-1,&\text{if $n\nmid j$,}
\end{cases}
$$
for any complex number $z$ with $|z|<1$, we have
\begin{align}
G(\zeta,z)=\sum_{k=1}^{n-1}\zeta^k \sum_{j=0}^\infty \zeta^{2jk}(j+1)z^j
&=\sum_{j=1}^{\infty}j\zeta^{j-1} \sum_{k=1}^{n-1} \zeta^{(2j-1)k}  \notag \\[5pt]
&=n\sum_{j=1}^\infty \left(jn-\frac{n-1}{2}\right) z^{jn-\frac{n+1}{2}}-\sum_{j=1}^\infty jz^{j-1} \notag \\[5pt]
&=\frac{n^2z^{\frac{n-1}{2}}}{(1-z^n)^2}-\frac{n(n-1)z^{\frac{n-1}{2}}}{2(1-z^n)}-\frac{1}{(1-z)^2}.  \label{eq:lomida}
\end{align}
Letting $z\to 1$ in \eqref{eq:lomida} and using L'ospital's rule, we get
\begin{align*}
G(\zeta)=G(\zeta,1)&=\lim_{z\to 1}\frac{2n^2z^{\frac{n-1}{2}}(1-z)^2-n(n-1)z^{\frac{n-1}{2}}(1-z^n)(1-z)^2-2(1-z^n)^2}{2(1-z^n)^2 (1-z)^2}\\[5pt]
&=\frac{n^2-1}{24},
\end{align*}
as desired.
\qed

\section{Proof of Theorem \ref{thm:q-Zudilin-0} }
\noindent{\it Proof of \eqref{eq:q-div-WZ-1}.} Define the following two functions in $q$:
\begin{align*}
F(n,k) &=[3n+2k+1]\frac{(q;q^2)_{n}(q^{2k+1};q^2)_{n}^2 q^{-{n+1\choose 2}-(2n+1)k} }{(q;q)_{n}^2 (q^2;q^2)_{n}}, \\[5pt]
G(n,k) &=-\frac{(1+q^{n+2k-1})(q;q^2)_{n}(q^{2k+1};q^2)_{n-1}^2 q^{-{n\choose 2}-(2n-1)k} }{(1-q)(q;q)_{n-1}^2 (q^2;q^2)_{n-1}},
\end{align*}
where we use the convention that $1/(q^2;q^2)_{a}=0$ for any negative integer $a$. Let $m$ be a positive odd integer. Then
\begin{align}
\sum_{k=0}^{\frac{m-1}{2}}[3k+1]\frac{(q;q^2)_k^3 q^{-{k+1\choose 2} } }{(q;q)_k^2 (q^2;q^2)_k}
=\sum_{n=0}^{\frac{m-1}{2}}F(n,0).  \label{eq:oursum}
\end{align}

It is easy to check that
\begin{align}
F(n,k-1)-F(n,k)=G(n+1,k)-G(n,k).  \label{eq:fnk-gnk}
\end{align}
Namely, the functions $F(n,k)$ and $G(n,k)$ form a $q$-WZ pair. (Usually the $q$-hypergeometric functions satisfying \eqref{eq:diff-fnk-gnk} are called a $q$-WZ pair.
But by symmetry it is convenient and reasonable to call the $q$-hypergeometric functions satisfying \eqref{eq:fnk-gnk} a $q$-WZ pair too, just as Zudilin \cite{Zudilin}).
Summing \eqref{eq:fnk-gnk} over $n=0,1,\ldots,\frac{m-1}{2}$, we obtain
\begin{align}
\sum_{n=0}^{\frac{m-1}{2}}F(n,k-1)-\sum_{n=0}^{\frac{m-1}{2}}F(n,k)=G\left(\frac{m+1}{2},k\right)-G(0,k)=G\left(\frac{m+1}{2},k\right).  \label{eq:fnk-gn0-00}
\end{align}
By \eqref{eq:mod-q2n-3}, for $k=1,2,\ldots, \frac{m-1}{2}$, we have
\begin{align}
G\left(\frac{m+1}{2},k\right) &=-\frac{(1+q^{\frac{m-1}{2}+2k})(q;q^2)_{(m+1)/2}(q^{2k+1};q^2)_{(m-1)/2}^2 q^{-\frac{m^2-1}{8}-pk} }{(1-q)(q;q)_{(m-1)/2}^3 (-q;q)_{(m-1)/2}} \notag  \\[5pt]
&\equiv 0\pmod{[m]\Phi_m(q)^2},  \label{eq:subs}
\end{align}
because $(-q;q)_{(m-1)/2}$ is relatively prime to $[m]$. Substituting \eqref{eq:subs} into \eqref{eq:fnk-gn0-00},
we see that
\begin{align*}
\sum_{n=0}^{\frac{m-1}{2}}F(n,0)\equiv \sum_{n=0}^{\frac{m-1}{2}}F(n,1) \equiv \sum_{n=0}^{\frac{m-1}{2}}F(n,2)\equiv\cdots
\equiv \sum_{n=0}^{\frac{m-1}{2}}F\left(n,\frac{m-1}{2}\right) \pmod{[m]\Phi_m(q)^2}.
\end{align*}
Hence, modulo $[m]\Phi_m(q)^2$, the sum \eqref{eq:oursum} can be replaced by
\begin{align}
\sum_{n=0}^{\frac{m-1}{2}}F\left(n,\frac{m-1}{2}\right)
&=\sum_{n=0}^{\frac{m-1}{2}}[3n+m]\frac{(q;q^2)_{n}(q^{m};q^2)_{n}^2 q^{-{n+1\choose 2}-\frac{(2n+1)(m-1)}{2}} }{(q;q)_{n}^2 (q^2;q^2)_{n}}  \notag \\[5pt]
&=[m]q^{\frac{1-m}{2}}+[m]\sum_{n=1}^{\frac{m-1}{2}}\frac{(q;q^2)_{n}(q^{m};q^2)_{n}^2 q^{-{n+1\choose 2}-\frac{(2n+1)(m-1)}{2}} }{(q;q)_{n}^2 (q^2;q^2)_{n}}  \notag \\[5pt]
&\quad+q^{m}\sum_{n=1}^{\frac{m-1}{2}}\frac{[3n](q;q^2)_{n}(q^{m};q^2)_{n}^2 q^{-{n+1\choose 2}-\frac{(2n+1)(m-1)}{2}} }{(q;q)_{n}^2 (q^2;q^2)_{n}}  \label{eq:reslut}
\end{align}

Similarly to the proof of Lemma \ref{eq:lem-2nn}, we can show that the denominator of the reduced form of  the fraction
$$
\frac{(q;q^2)_{n}(q^{m};q^2)_{n}^2 }{(q;q)_{n}^2 (q^2;q^2)_{n}}
$$
is relatively prime to $[m]$ (since $m$ is odd). It is clear that the numerator of this reduced
form is divisible by $\Phi_m(q)^2$ for $1\leqslant n\leqslant m-1$.
Thus, comparing the expression \eqref{eq:reslut} with \eqref{eq:oursum}, we see that \eqref{eq:q-div-WZ-1} with $n=m$
is equivalent to
\begin{align}
q^{m}\sum_{n=1}^{\frac{m-1}{2}}\frac{[3n](q;q^2)_{n}(q^{m};q^2)_{n}^2 q^{-{n+1\choose 2}-\frac{(2n+1)(m-1)}{2}} }{(q;q)_{n}^2 (q^2;q^2)_{n}}
\equiv 0\pmod{[m]\Phi_m(q)^2}.  \label{eq:reduce-3}
\end{align}
We can further show that
\begin{align*}
\frac{[3n](q;q^2)_{n}(q^{m};q^2)_{n}^2 }{(q;q)_{n}^2 (q^2;q^2)_{n}}\equiv 0 \pmod{[m]} \quad\text{for $n\leqslant\frac{m-1}{2}$}.  
\end{align*}
Thus, writing the left-hand side of \eqref{eq:reduce-3} as
\begin{align*}
[m]^2\sum_{n=1}^{\frac{m-1}{2}}[3n]\frac{(1-q)^2(q;q^2)_{n}(q^{m+2};q^2)_{n-1}^2 q^{m-{n+1\choose 2}-\frac{(2n+1)(m-1)}{2}} }{(q;q)_{n}^2 (q^2;q^2)_{n}},
\end{align*}
and noticing that $[m]\equiv 0\pmod{\Phi_m(q)}$, we see that \eqref{eq:reduce-3} is equivalent to
\begin{align}
\sum_{n=1}^{\frac{m-1}{2}}[3n]\frac{(1-q)^2(q;q^2)_{n}(q^{m+2};q^2)_{n-1}^2 q^{m-{n+1\choose 2}-\frac{(2n+1)(m-1)}{2}} }{(q;q)_{n}^2 (q^2;q^2)_{n}}
\equiv 0\pmod{\Phi_m(q)}. \label{eq:reduce}
\end{align}
Since $q^m\equiv 1\pmod{\Phi_m(q)}$, we can reduce \eqref{eq:reduce} to its equivalent form
\begin{align*}
&\sum_{n=1}^{\frac{m-1}{2}}(1+q^n+q^{2n})\frac{(1-q)(q;q^2)_{n}(q^{2};q^2)_{n-1} q^{-{n\choose 2}-\frac{m-1}{2}} }{(1+q^n)(q;q)_{n}^2 } \\[5pt]
&\quad =\sum_{n=1}^{\frac{m-1}{2}}\frac{1+q^n+q^{2n}}{(1+q^n)[2n]}{2n\brack n}q^{-{n\choose 2}-\frac{m-1}{2} }
\equiv 0\pmod{\Phi_m(q)},
\end{align*}
which is the $n=m$ case of \eqref{eq:3k2kk} differing only by a factor $q^{-\frac{m-1}{2}}$.
\qed

\medskip
\noindent{\it Proof of \eqref{eq:q-div-WZ-2}.} Let $m$ be a positive odd integer again. Then
\begin{align}
\sum_{k=0}^{m-1}[3k+1]\frac{(q;q^2)_k^3 q^{-{k+1\choose 2} } }{(q;q)_k^2 (q^2;q^2)_k}
=\sum_{n=0}^{m-1}F(n,0).  \label{eq:oursum-2}
\end{align}
Similarly as before, summing \eqref{eq:fnk-gnk} over $n=0,1,\ldots,m-1$, we obtain
\begin{align}
\sum_{n=0}^{m-1}F(n,k-1)-\sum_{n=0}^{m-1}F(n,k)=G(m,k).  \label{eq:fnk-gn0-66}
\end{align}
By \eqref{eq:mod-q2n}, for $k=1,2,\ldots, m-1$, we have
\begin{align*}
G(m,k) &=-\frac{(1+q^{m+2k-1})(q;q^2)_{m}(q^{2k+1};q^2)_{m-1}^2 q^{-{m\choose 2}-(2m-1)k} }{(1-q)(q;q)_{m-1}^3 (-q;q)_{m-1}}
\equiv 0\pmod{[m]\Phi_m(q)^2},
\end{align*}
since $\gcd((-q;q)_{m-1},[m])=1$. Thus, we get
\begin{align*}
\sum_{n=0}^{m-1}F(n,0)\equiv \sum_{n=0}^{m-1}F(n,1) \equiv\cdots\equiv \sum_{n=0}^{m-1}F\left(n,\frac{m-1}{2}\right) \pmod{[m]\Phi_m(q)^2}.
\end{align*}
Therefore, modulo $[m]\Phi_m(q)^2$, the sum \eqref{eq:oursum-2} can be replaced by
\begin{align*}
\sum_{n=0}^{m-1}F\left(n,\frac{m-1}{2}\right)
=\sum_{n=0}^{m-1}[3n+m]\frac{(q;q^2)_{n}(q^{m};q^2)_{n}^2 q^{-{n+1\choose 2}-\frac{(2n+1)(m-1)}{2}} }{(q;q)_{n}^2 (q^2;q^2)_{n}}.
\end{align*}
Similarly to the proof of \eqref{eq:q-div-WZ-1},  we see that \eqref{eq:q-div-WZ-2} with $n=m$
is equivalent to
\begin{align*}
\sum_{n=1}^{m-1}[3n]\frac{(1-q)^2(q;q^2)_{n}(q^{m+2};q^2)_{n-1}^2 q^{m-{n+1\choose 2}-\frac{(2n+1)(m-1)}{2}} }{(q;q)_{n}^2 (q^2;q^2)_{n}}
\equiv 0\pmod{\Phi_m(q)},
\end{align*}
which immediately follows again from the $n=m$ case of \eqref{eq:3k2kk} (note that we may sum in \eqref{eq:3k2kk} up to $n-1$).
\qed

\medskip
\noindent{\it Remark.} The functions $F(n,k)$ and $G(n,k)$ are not easy to find (they always requires a preliminary human
guess),
but once they are guessed out correctly they can hopefully be
proved via the $q$-WZ method (see \cite{Mohammed, Zeilberger}). Since Guillera and Zudilin \cite{GuZu} have given the corresponding WZ pair (the $q=1$ case),
we can find the $q$-WZ pair not so difficultly.

\section{Proof of Theorem \ref{thm:q-Zudilin-3}}
Let
\begin{align*}
F(n,k) &=(-1)^{n}[3n-2k+1]{2n-2k\brack n}\frac{(q;q^2)_{n}(q;q^2)_{n-k} }{(q;q)_{n} (q^2;q^2)_{n-k}}, \\[5pt]
G(n,k) &=(-1)^{n+1}[n]{2n-2k\brack n-1}\frac{(q;q^2)_{n}(q;q^2)_{n-k} q^{n+1-2k} }{(q;q)_{n} (q^2;q^2)_{n-k}}.
\end{align*}
These functions $F(n,k)$ and $G(n,k)$ are $q$-analogues of the WZ pair given by He in the proof of \cite[Theorem 1.1]{He}.
It is easy to check that $F(n,k)$ and $G(n,k)$ satisfy
\begin{align}
F(n,k-1)-F(n,k)=G(n+1,k)-G(n,k).  \label{eq:fnk-gnk-new}
\end{align}
Let $m$ be a positive odd integer. Since
$$
{2k\brack k}
=\frac{(q;q^2)_k (-q;q)_k^2}{(q^2;q^2)_k},
$$
we have
\begin{align}
\sum_{k=0}^{m-1}(-1)^k [3k+1]\frac{(q;q^2)_k^3}{(q;q)_k^3}
=\sum_{n=0}^{m-1}F(n,0).  \label{eq:fn0}
\end{align}

Summing \eqref{eq:fnk-gnk-new} over $n=0,1,\ldots,m-1$, we obtain
\begin{align}
\sum_{n=0}^{m-1}F(n,k-1)-\sum_{n=0}^{m-1}F(n,k)=G(m,k).  \label{eq:fnk-gn0-00}
\end{align}
By Lemma \ref{lem:big-qbino}, for $k=1,\ldots,\frac{m-1}{2}$, we have $G(m,k)\equiv 0\pmod{[m]\Phi_m(q)^2}$. Therefore, from \eqref{eq:fnk-gn0-00} we deduce that
\begin{align}
\sum_{n=0}^{m-1}F(n,0)\equiv \sum_{n=0}^{m-1}F(n,1) \equiv\cdots\equiv \sum_{n=0}^{m-1}F\left(n,\frac{m-1}{2}\right) \pmod{[m]\Phi_m(q)^2}.  \label{eq:fn1-fn2}
\end{align}
Furthermore, since $F(n,k)=0$ for $n<2k$, we get
\begin{align}
\sum_{n=0}^{m-1}F\left(n,\frac{m-1}{2}\right)=F\left(m-1,\frac{m-1}{2}\right)
&=[2m-1]\frac{(q;q^2)_{m-1}(q;q^2)_{(m-1)/2} }{(q;q)_{m-1} (q^2;q^2)_{(m-1)/2}}  \notag\\[5pt]
&={2m-1\brack m-1}{m-1\brack \frac{m-1}{2}}_{q^2}\frac{[m]}{(-q;q)_{m-1}^2}.  \label{eq:newhe}
\end{align}
Substituting the following two congruences (see \cite[(3.1)]{Guo2018} and \cite[(1.5)]{LPZ})
\begin{align*}
{2m-1\brack m-1}
&\equiv (-1)^{m-1}q^{m\choose 2} \pmod{\Phi_m(q)^2}, \\[5pt]
{m-1\brack \frac{m-1}{2}}_{q^2}
&\equiv (-1)^{\frac{m-1}{2}}q^{\frac{1-m^2}{4}}(-q;q)_{m-1}^2  \pmod{\Phi_m(q)^2}
\end{align*}
into \eqref{eq:newhe}, and noticing \eqref{eq:fn0} and \eqref{eq:fn1-fn2}, we complete the proof of the theorem for $n=m$.

\medskip\noindent{\it Remark.} The $q=1$ case of our proof gives a new proof of \eqref{eq:div-3}, which is simpler than Guillera and Zudilin's original proof.
Moreover, summing \eqref{eq:fnk-gnk-new} over $n=0,1,\ldots,\frac{m-1}{2}$, we can only show that
\begin{align}
\sum_{k=0}^{\frac{m-1}{2}}(-1)^k [3k+1]\frac{(q;q^2)_k^3}{(q;q)_k^3}
\equiv 0 \pmod{\Phi_m(q)}.  \label{eq:remark-Zudilin}
\end{align}
Therefore, to confirm Conjecture \ref{conj:q-Zudilin-3} we need new methods or techniques, though Theorem~\ref{thm:q-Zudilin-3} and Conjecture~\ref{conj:q-Zudilin-3}
are equivalent when $n$ is an odd prime power.

\section{Proof of Theorem \ref{thm:q-Zudilin-MZ}}
Summing \eqref{eq:fnk-gnk-new} over $n$ from $0$ to $N$, we obtain
\begin{align}
\sum_{n=0}^{N}F(n,k-1)-\sum_{n=0}^{N}F(n,k)=G\left(N+1,k\right).  \label{eq:fnk-nnn}
\end{align}
Furthermore, summing \eqref{eq:fnk-nnn} over $k$ from $1$ to $N$, we get
\begin{align}
\sum_{n=0}^{N}F(n,0)-\sum_{n=0}^{N}F(n,N)=\sum_{k=1}^{N}G\left(N+1,k\right).  \label{eq:sum-fn0}
\end{align}
Since  $F(n,N)=0$ for $n=0,1,\ldots,N$, it follows from \eqref{eq:sum-fn0} that
\begin{align}
&\hskip -3mm (-q;q)_N^3\sum_{n=0}^{N}F(n,0) \notag\\[5pt]
&=(-q;q)_N^3\sum_{k=1}^{N}G\left(N+1,k\right)  \notag\\[5pt]
&=(-1)^N\sum_{k=1}^{N}[N+1]{2N+2\brack N+1}{2N-2k+2\brack N}{2N-2k+2\brack N-k+1}\frac{(-q;q)_N^2 q^{N-2k+2}}{(1+q^{N+1})(-q;q)_{N-k+1}^2}.  \label{eq:sun-qbino}
\end{align}
Since $[N+1]{2N+2\brack N+1}/(1+q^{N+1})=[2N+1]{2N\brack N}$, it is easy to see that each summand on the right-hand side of \eqref{eq:sun-qbino}
is divisible by $(1+q^N)^2[2N+1]{2N\brack N}$ (whenever $k=1$ or $k\geqslant 2$). This proves that the congruence \eqref{eq:q-Sun-2} holds for $n=N$.  \qed

\section{Proof of Theorem \ref{thm:q-st}}
Replacing $n$ by $n-1$ in \eqref{eq:q-Staver}, we obtain
\begin{align}
\sum_{k=1}^{n-1}\frac{[3k]}{[2k]^2}{2k\brack k}q^{-{k\choose 2}}
=[n]{2n-1\brack n-1}\sum_{k=1}^{n-1}\frac{q^{-{n-2k\choose 2}}}{[2k]^2{n-1\brack k}_{q^2}^2}. \label{eq:qst-n-1}
\end{align}
It is easy to see that ${2n-1\brack n-1}\equiv 1\pmod{\Phi_n(q)}$ and ${n-1\brack k}_{q^2}\equiv (-1)^k q^{-k^2-k}\pmod{\Phi_n(q)}$.
Since $[n]\equiv 0\pmod{\Phi_n(q)}$ and $q^n\equiv 1\pmod{\Phi_n(q)}$, we deduce from \eqref{eq:qst-n-1} that
\begin{align*}
\sum_{k=1}^{n-1}\frac{[3k]}{[2k]^2}{2k\brack k}q^{-{k\choose 2}}
\equiv [n]\sum_{k=1}^{n-1}\frac{q^{k}}{[2k]^2} \pmod{\Phi_n(q)^2}.
\end{align*}
The proof then follows from \eqref{eq:sp-sim}.

\section{Concluding remarks and open problems}
Numerical calculation suggests the following refinement of \eqref{eq:q-div-WZ-2}.
\begin{conj}\label{conj:q-Zudilin-0}
Let $n$ be a positive odd integer. Then
\begin{align*}
\sum_{k=0}^{n-1}[3k+1]\frac{(q;q^2)_k^3 q^{-{k+1\choose 2} } }{(q;q)_k^2 (q^2;q^2)_k}
&\equiv [n]q^{\frac{1-n}{2}}+\frac{(n^2-1)(1-q)^2}{24}[n]^3 q^{\frac{1-n}{2}} \pmod{[n]\Phi_n(q)^3}.
\end{align*}
\end{conj}
Note that if $n=p^r$ and $q=1$, then the above congruence reduces to
\begin{align}
\sum_{k=0}^{p^r-1} \frac{(\frac{1}{2})_k^3}{k!^3}(3k+1)2^{2k} &\equiv p^r\pmod{p^{r+3}}\quad\text{for $p>3$}. \label{eq:sun-hu}
\end{align}
A stronger version of \eqref{eq:sun-hu} modulo $p^{r+4}$ was conjectured by Sun\cite[Conjecture 5.1(ii)]{Sun4}. Moreover, the $r=1$ case of \eqref{eq:sun-hu} has been proved by Dian-Wang Hu in his Ph.D. thesis.

Motivated by \eqref{eq:q-Hamme}, we give  the following $q$-analogue of a supercongruence in \cite[Conjecture 5.9]{Sun4} (with $a=1$).
\begin{conj}Let $n$ be a positive integer with $n\equiv 1\pmod 4$. Then
\begin{align*}
\sum_{k=0}^{\frac{n-1}{2}}[4k+1]\frac{(q;q^2)_k^3}{(q^2;q^2)_k^3} q^{\frac{k(n^2-2nk-n-2)}{4}}\equiv 0\pmod{\Phi_n(q)^2}.
\end{align*}
\end{conj}

 Guillera and Zudilin \cite{GuZu} proved \eqref{eq:div-3} by using the following WZ pair:
\begin{align*}
F(n,k)&=(-1)^n(3n+2k+1)\frac{(\frac{1}{2})_n(\frac{1}{2}+k)_n^2(\frac{1}{2})_k }{(1)_n^2(1+2k)_n(1)_k}2^{3n},\\[5pt]
G(n,k)&=(-1)^n \frac{(\frac{1}{2})_n(\frac{1}{2}+k)_{n-1}^2(\frac{1}{2})_k }{(1)_{n-1}^2(1+2k)_{n-1}(1)_k}2^{3n-2}.
\end{align*}
We also found a $q$-analogue of the above WZ pair. Unfortunately, we cannot find and prove the corresponding $q$-analogue of \cite[Lemma 4]{GuZu}, and so we are unable to
use this $q$-WZ pair to prove Theorem \ref{thm:q-Zudilin-3}. On the other hand, from Theorem \ref{thm:q-Zudilin-3} and this $q$-WZ pair we can deduce that,
for any positive odd integer $n$, there holds
\begin{align*}
\sum_{k=1}^{n-1}(-1)^k\frac{[n+3k](-q^{n+1};q)_{k-1}}{[2k]}{2k\brack k}
\equiv (1+q^n)\left(q^{n\choose 2}-(-q;q)_{n-1}\right)   \pmod{\Phi_n(q)^2},
\end{align*}
which seems difficult to be proved directly. Moreover, using this $q$-WZ pair, we can show that the congruence \eqref{eq:conj-Zudilin-3} modulo $[n]$ is true, which is
a little stronger than \eqref{eq:remark-Zudilin}. We plan to give details of the proofs of these two congruences in another paper.

Mao and Zhang \cite{MZ} proved \eqref{eq:sun-bino-1} by applying the $x=-\frac{1}{2}$ case of the following identity (see \cite[Lemma 3.2]{MS}):
\begin{align}
\sum_{k=0}^{n}{n\choose k}^2{x+k\choose 2n+1}
=\frac{1}{(4n+2){2n\choose n}}\sum_{k=0}^{n}(2x-3k){x\choose k}^2{2k\choose k}.  \label{eq:Mao-Sun}
\end{align}
Therefore, a possible way to prove Conjecture \ref{conj:q-Zudilin-MZ} is to give a $q$-analogue of the curious binomial coefficient identity \eqref{eq:Mao-Sun}.

Swisher \cite{Swisher} has proposed many interesting conjectures on supercongruences that generalize van Hamme's original 13 supercongruences \cite{Hamme}.
For example, she conjectured that \cite[(J.3)]{Swisher}
\begin{align*}
\sum_{k=0}^{\frac{p^r-1}{2}} \frac{(6k+1)(\frac{1}{2})_k^3}{k!^3 4^k}
\equiv (-1)^\frac{p-1}{2}p \sum_{k=0}^{\frac{p^{r-1}-1}{2}} \frac{(6k+1)(\frac{1}{2})_k^3}{k!^3 4^k} \pmod{p^{4r}} \quad\text{for $p>3$}.
\end{align*}

Motivated by Swisher's conjectures and also the conjectures of Sun \cite{Sun}, we propose the following refinements of
\eqref{eq:div-gen-r-1}, \eqref{eq:div-gen-r-2}, \eqref{eq:div-3-pr}, and \eqref{eq:sun-hu} (there is also a similar conjecture in \cite{OZ}).

\begin{conj}Let $p$ be an odd prime and $r$ a positive integer. Then
\begin{align*}
\sum_{k=0}^{\frac{p^r-1}{2}} \frac{(\frac{1}{2})_k^3}{k!^3}(3k+1)2^{2k}
&\equiv p \sum_{k=0}^{\frac{p^{r-1}-1}{2}} \frac{(\frac{1}{2})_k^3}{k!^3}(3k+1)2^{2k} \pmod{p^{3r}}, \\[5pt]
\sum_{k=0}^{p^r-1} \frac{(\frac{1}{2})_k^3}{k!^3}(3k+1)2^{2k}
&\equiv p \sum_{k=0}^{p^{r-1}-1} \frac{(\frac{1}{2})_k^3}{k!^3}(3k+1)2^{2k} \pmod{p^{4r-\delta_{p,3}}}, \\[5pt]
\sum_{k=0}^{\frac{p^r-1}{2}} \frac{(\frac{1}{2})_k^3}{k!^3}(3k+1)(-1)^k 2^{3k}
&\equiv p(-1)^{\frac{p-1}{2}} \sum_{k=0}^{\frac{p^{r-1}-1}{2}} \frac{(\frac{1}{2})_k^3}{k!^3}(3k+1)2^{2k} \pmod{p^{3r+\delta_{p,3}}}, \\[5pt]
\sum_{k=0}^{p^r-1} \frac{(\frac{1}{2})_k^3}{k!^3}(3k+1)(-1)^k 2^{3k}
&\equiv p(-1)^{\frac{p-1}{2}} \sum_{k=0}^{p^{r-1}-1} \frac{(\frac{1}{2})_k^3}{k!^3}(3k+1)2^{2k} \pmod{p^{3r}},
\end{align*}
where $\delta$ is the Kronecker delta with $\delta_{i,j}=1$ if $i=j$ and $\delta_{i,j}=0$ otherwise.
\end{conj}

Finally, it would be interesting to find a $q$-analogue of \eqref{eq:div-2}.

\vskip 2mm \noindent{\bf Acknowledgments.} This work was partially
supported by the National Natural Science Foundation of China (grant 11771175),
the Natural Science Foundation of Jiangsu Province (grant BK20161304),
and the Qing Lan Project of Education Committee of Jiangsu Province.

\end{document}